\documentclass{amsart}
\usepackage{amsfonts}
\usepackage{amscd}
\usepackage{amsbsy}
\usepackage{amssymb}
\usepackage{graphicx}
\usepackage[all]{xy}

\newtheorem{teo}{Theorem}[section]

\newtheorem{cor}[teo]{Corollary}
\newtheorem{lem}[teo]{Lemma}
\newtheorem{rem}[teo]{Remark}
\newtheorem{defn}[teo]{Definition}

\newtheorem{quest}[teo]{Question}

\newtheoremstyle{example}{3pt}{3pt}%
{}
{}
{\sc}
{:}
{.5em}
{}
\theoremstyle{example}

\def\CC{{\mathcal C}}
\def\OO{{\mathcal O}}
\def\UU{{\mathcal U}}
\def\DD{{\mathcal D}}

\def\LL{{\mathcal L}}

\def\OO{{\mathcal O}}

\def\RR{{\mathcal R}}
\def\SS{{\mathcal S}}

\def\UU{{\mathcal U}}
\def\VV{{\mathcal V}}

\def\ZZ{{\mathcal Z}}

\def\fract#1#2{\raise4pt\hbox{$ #1 \atop #2 $}}

\def\C{{\mathbb C}}

\def\H{{\mathbb H}}

\def\O{{\mathbb O}}

\def\R{{\mathbb R}}

\def\Z{{\mathbb Z}}

\def\gg{{\mathfrak g}}
\def\gh{{\mathfrak h}}

\def\gm{{\mathfrak m}}

\def\gs{{\mathfrak s}}
\def\gt{{\mathfrak t}}
\def\gu{{\mathfrak u}}

\def\gH{{\mathfrak H}}

\def\la#1{\hbox to #1pc{\leftarrowfill}}
\def\ra#1{\hbox to #1pc{\rightarrowfill}}

\def\Qquot{\setminus\!\!\setminus\!\!\setminus}

\begin{document}
\bibliographystyle{amsalpha}
\title{Quaternionic K\"ahler Reductions of Wolf Spaces}
\author{Daniele Grandini}
\address{Department of Mathematics \& Statistics,
University of New Mexico,
Albuquerque, NM 87131.}
\address{Dipartimento di Matematica "F.Enriques",
Universit\`a degli Studi di Milano, Milano, Italy.}
\email{grandini@mat.unimi.it, daniele@math.unm.edu}
\maketitle

\noindent \emph{Abstract:} The main purpose of the following
article is to introduce a \emph{Lie theoretical} approach to the
problem of classifying pseudo quaternionic-K\"ahler (QK)
reductions of the pseudo QK symmetric spaces, otherwise
called \emph{generalized Wolf spaces}.\\

The history of QK geometry starts with the celebrated
\emph{Berger's Theorem} \cite{Ber55} which classifies all the
irreducible holonomy groups for not locally symmetric
pseudo-riemannian manifolds. In fact, a pseudo QK manifold $(M,g)$
of dimension $4n$ ($n>1$) is traditionally defined by the
reduction of the holonomy group to a subgroup of ${\rm
Sp}(k,l){\rm Sp}(1)$ ($k+l=n$). Alekseevsky proved \cite{Al2} that
any pseudo QK manifold is Einstein and satisfies some additional
curvature condition, so it is possible to extend in a natural way
the definition of QK manifolds to 4-manifolds: an oriented
4-manifold is said to be QK if it is Einstein and self-dual.
Furthermore, the whole definition can be naturally extended to
orbifolds \cite{GaLa88}.\\
Examples of pseudo QK manifolds are not too many, and most of them
are homogeneous spaces. In particular, Alekseevsky proved that all
homogeneous, Riemannian QK manifolds with positive scalar
curvature are (compact and) symmetric \cite{Al2}. These spaces
were classified by Wolf \cite{Wol65} and have been called
\emph{Wolf spaces}.\\
The Wolf spaces together with their duals can be characterized as
all the non-scalar flat QK manifolds admitting a transitive
unimodular group of isometries \cite{AlCor97}.\\
The duals of Wolf spaces don't classify all the homogeneous,
Riemannian QK manifolds with negative scalar curvature: in fact,
there are more such spaces, like the so-called \emph{Alekseevskian
spaces} \cite{Al2}.\\Most examples of non-homogeneous  QK
\emph{orbifolds} emerge via the so-called \textbf{symmetry
reduction}. This process can be seen as a variation of a
well-known construction of Marsden and Weinstein developed in the
context of Poisson and symplectic geometry (see \cite{MaWe74} or
\cite{AbMa78}). The Marsden-Weinstein quotient considers a
symplectic manifold with some symmetries. A \emph{new} symplectic
manifold of lower dimension and fewer symmetries is then obtained
by ``dividing out" some symmetries in a ``symplectic fashion".
This simple idea has been more recently applied in many different
geometric situations. Various generalizations of the symplectic
reductions include K\"ahler quotients, hyperk\"ahler and
hypercomplex quotients, quaternionic K\"ahler and quaternionic
quotients, 3-Sasakian, Sasakian and contact quotients to mention
just a few. The so-called \textbf{QK reduction} has been
introduced by Galicki and Lawson (see \cite{Ga86, GaLa88}). Here
one starts with a QK space $M^{4n}$ with some symmetry $H$. A new
QK space of dimension $4n-4{\rm dim}(H)$ is constructed
out of the quaternionic K\"ahler moment map.\\

In this paper we specialize the QK-reduction to pseudo QK
symmetric spaces, that we call \emph{generalized Wolf spaces} (or
sometimes just \emph{Wolf spaces}). Recently, this spaces have
been classified by Alekseevsky and Cort{\'e}s \cite{AlCor05}.\\
In particular, we show that given a (generalized) Wolf space $G/H$
the quaternionic K\"ahler moment map can be lifted to a
$\R^3$-valued map defined on the group $G$.

In the first two sections of the paper we briefly introduce the
pseudo QK geometry and the fundamental notion of \emph{diamond
diagram} $\diamondsuit$,  a ``bundle diagram" built up by Boyer,
Galicki and Mann \cite{BGM93a} which functorially relates pseudo
QK structures to pseudo 3-Sasakian, pseudo hyperk\"ahler and
pseudo twistor structures, studied respectively by Konishi
\cite{Kon74}, Swann
\cite{Sw91} and Salamon \cite{Sal82}.\\
In the third section we sketch the QK reduction and show that the
QK moment map is canonically associated to a moment map of the
whole diamond diagram \cite{BGM93a, Sw91}.\\
In the fourth section
we introduce the Wolf spaces and describe
their associate diamond diagrams.\\
In the fifth and sixth sections we give the explicit formula for
the ``diamond moment map" associated to Wolf spaces, though as a
$\R^3$-valued map defined on the total group $G$. We discuss its
properties and focus our attention to 1-dimensional actions:\\in
this case, we show that the classification of reductions is
related to the problem to classify the adjoint orbits of the total
group.\\
In the last part of the paper we  apply the whole construction to
the Grassmann manifold ${\rm SO}(7)/{\rm SO}(3)\times {\rm
SO}(4)$, together its dual space ${\rm SO}_0(3,4)/{\rm
SO}(3)\times {\rm SO}(4)$.\\
Firstly, we find  ``canonical" reductions of them, related to some
\emph{normed algebra} structures of $\R^8$, i.e. the reductions
$$\xymatrix{\diamondsuit({\rm SO}(7)/{\rm SO}(3)\times {\rm
SO}(4))\ar @{=>}[r]^{\ \ \ \ \ \ {\rm
U}(1)}&\diamondsuit(\Z_3\setminus G_2/{\rm SO}(4))}$$
$$\xymatrix{\diamondsuit({\rm SO}_0(3,4)/{\rm SO}(3)\times {\rm
SO}(4))\ar @{=>}[r]^{\ \ \ \ \ \ {\rm
U}(1)}&\diamondsuit(\Z_3\setminus G_{2(2)}/{\rm SO}(4))}$$ Next,
we classify the adjoint orbits of ${\rm SO}(7)$ and ${\rm
SO}(3,4)$.\\
Since the compactness of ${\rm SO}(7)$, the classification of its
adjoint orbits is elementary and the reductions are nothing but
the ``weighted deformations" of the canonical reduction.\\
In the case of ${\rm SO}(3,4)$, the classification of the adjoint
orbits is much more complicated and, besides the weighted
deformations of the canonical reduction, we obtain many other
reductions which are always \emph{smooth}.

\section{Introduction}
\noindent Let $\OO^{4n}$ denote a $4n$-dimensional (connected)
manifold (or orbifold). An \textbf{almost quaternionic structure}
on $\OO^{4n}$ is a rank 3 bundle ${\gH}\subset {\rm End}(T\OO)$
such that at every point there are local sections $J_1,J_2,J_3$ of
$\gH$ satisfying
$$J_aJ_b=-\delta_{ab}+\epsilon^{abc}J_c.$$
It follows that  $\OO^{4n}$ is almost quaternionic if and only if
the structure group of the tangent bundle reduces to ${\rm
GL}(n,\H)\H^*={\rm GL}(n,\H){\rm Sp}(1).$ Given a pseudoriemannian
metric $g$ on $\OO^{4n}$, the pair $(\gH,g)$ is called
\textbf{pseudo-hyperhermitian structure} on $\OO^{4n}$ if the
bundle $\gH$ defines an almost quaternionic structure as defined
above and $\gH$ is contained in the vector subbundle of ${\rm
End}(T\OO)$ generated by the isometries of $g$:
$${\gH}\subset\R\cdot{\rm O}_g(T\OO).$$ In terms of the local frame
$J_1,J_2,J_3$, this means that
$$g(J_aX,J_aY)=g(X,Y)$$
for any local vector fields $X,Y$ and for all $a=1,2,3$. In
particular, the
signature of $g$ has to be of the type $(4k,4l)$, $k+l=n$.\\
\begin{rem}Given a pseudo-hyperhermitian structure $(\gH,g)$, one can use the local duality
$\xy\xymatrix{J_a \ar@{<->}[r]^g& \omega_a}\endxy$ to identify
$\gH$ with a subbundle of $ \bigwedge^2T^*M$. In particular, the
local forms $\omega_a$ are almost symplectic.
\end{rem}
\noindent If $n>1$, $\OO^{4n}$ is said to be
\textbf{quaternionic-K\"ahler (QK)} with respect to the
pseudo-hyperhermitian structure $(\gH,g)$ if the the bundle $\gH$
is preserved by the Levi-Civita connection of $g$.\\
\begin{teo}Let $\OO^{4n}$ as above ($n\geq 1$) and let $g$ be a pseudoriemannian metric on $\OO$ with signature $(4k,4l)$
Then, there exists an almost quaternionic structure $\gH$ on $\OO$
such that $(\gH,g)$ is a pseudo-hyperhermitian structure if and
only if the tangent bundle reduces to the group ${\rm Sp}(k,l){\rm
Sp}(1)$. Moreover, if $n>1$, $\OO^{4n}$ is QK with respect to
$(\gH,g)$ if and only if $${\rm Hol}(g)\subset {\rm Sp}(k,l){\rm
Sp}(1)$$\end{teo} \noindent Given a pseudo-hyperhermitian
structure $(\gH,g)$ on $\OO^{4n}$, we can define two \emph{global}
tensors, locally written as
\begin{eqnarray}\Theta:=\sum_a\omega_a\otimes
J_a&\in&\displaystyle\Gamma\left(\wedge^2T^*\OO\otimes\gH\right),\\
\Omega:=\sum_a\omega_a\wedge\omega_a&\in&\displaystyle\Gamma\left(\wedge^4T^*\OO\right).\end{eqnarray}
where $\omega_a$ are local almost symplectic forms, dual to $J_a$
via the metric $g$.
\begin{teo} Let $(\gH,g)$ be a pseudo-hyperhermitian
structure on $\OO^{4n}$. If $n>1$, then the following conditions
are equivalent:
\begin{enumerate} \item $\OO^{4n}$ is QK w.r.t. $(\gH,g)$.%
\item $\Theta$ is parallel. %
\item $\Omega$ is parallel.
\end{enumerate}
\end{teo}
\noindent It follows that the bundle $\gH$ is orientable and the
pseudoriemannian metric $g$ induces a metric on $\gH$, given by
$$<J,J'>:=-\frac{1}{4n}{\rm trace}_g(J\circ J')$$
which is clearly positive defined. We can introduce the bundle of
the oriented orthogonal frames of $\gH$
$$\pi:\SS\rightarrow\OO,$$
called the \textbf{Konishi bundle} \cite{Kon74} of the
pseudo-hyperhermitian
structure $(\gH,g)$.\\
If $n>1$ and $\OO^{4n}$ is QK relative to $(\gH,g)$, then the
Levi-Civita connection induced on $\gH$ is metric relative to
$<,>$, which means that
$$X<J,J'>=<\nabla_XJ,J'>+<J,\nabla_XJ'>.$$
Moreover, the local coefficients arising from the formula
$$\nabla J_i=\alpha^{j}_i\otimes J_j$$
(or, equivalently, $\nabla \omega_i=\alpha^{j}_i\omega_j$ by
duality) define a connection on the ${\rm SO}(3)$-principal bundle
$\SS$, whose curvature $F$ is given by
$$F_i^j:=d\alpha_i^j-\alpha_i^k\wedge\alpha_k^j.$$
The curvature $F$ is related to the Riemann curvature $R$, seen as
a map
$$R:\bigwedge^2T\OO\rightarrow {\rm SkewEnd}(T\OO).$$
We get
$$[R,J_i]=F_i^j\otimes J_j.$$
As the image of $R$ takes values in
$\mathfrak{sp}(k,l)\oplus\mathfrak{sp}(1)$, we must have
$$[R,J_i]=[R_{\mathfrak{sp}(1)},J_i],$$
and, hence,
$$R_{\mathfrak{sp}(1)}=\frac{1}{2}(F_2^3\otimes J_1+F_3^1\otimes J_2+F_1^2\otimes J_3).$$
Furthermore, it turns out that
$$F_i^j=\lambda\epsilon^{ijk}\omega_k$$
for some constant $\lambda$ which gives

\begin{equation}\label{1}
R_{\mathfrak{sp}(1)}=\frac{\lambda}{2}\Theta\end{equation}
Finally, the condition $[R_{\mathfrak{sp}(n)},J_i]=0$ implies
$${\rm Ric}_g(X,Y)={\rm trace}_g\{u\mapsto
R_{\mathfrak{sp}(1)}(u,X)Y\}=-\frac{3\lambda}{2}g(X,Y).$$ In
particular, the metric $g$ is Einstein (see also \cite{Al1} or \cite{Sal82}).\\
Let us now discuss the $4$-dimensional case: clearly, $\OO^4$
admits a pseudo-hyperhermitian structure if and only if it is
oriented. In this case, if $(g,\gH)$ is a pseudo-hyperhermitian
structure on $\OO^4$, then $\pm g$ is a Riemannian metric and
$\gH$ is canonically identified with the bundle $\Lambda_+$ of the
self-dual 2-forms. Moreover, the holonomy is contained in ${\rm
SO}(4)={\rm Sp}(1){\rm Sp}(1)$, and the identities
${\nabla}\Omega={\nabla}\Theta=0$ are trivially satisfied.
\begin{defn}Let $(g,\Lambda_+)$ be a pseudo-hyperhermitian
structure on $\OO^4$. Then, $\OO^4$ is QK with respect to
$(g,\Lambda_+)$ if $g$ is Einstein and the condition (\ref{1}) is
satisfied.
\end{defn}
\noindent In particular, a QK-manifold has constant scalar
curvature and the scalar-flat QK manifolds are \emph{locally}
hyperk\"ahler.\\

 From this point on $\OO$ will always be a
non scalar-flat QK space of dimension greater or equal than $4$
and signature $(4k,4l)$.

\section{The diamond diagram}
\subsection{The $3$-Sasakian structure of the Konishi bundle}
\noindent Let $\eta=(\eta_1,\eta_2,\eta_3)$ be the connection
1-form of the Konishi bundle $\pi:\SS\rightarrow\OO$, defined
locally by the forms $\alpha_i^j$. It can be shown that the forms
$\eta_i$ are indeed \emph{contact forms} and they define uniquely
three global vector fields $\xi_1,\xi_2,\xi_3$ such that
$$\eta_i(\xi_i)\equiv 1,\qquad i_{\xi_i}d\eta_i=0.$$
Furthermore, one sees that the $\xi_i$'s define a global 3-frame
on $\SS$ (called \textbf{Reeb $3$-distribution}), satisfying the
rule
$$[\xi_a,\xi_b]=2\epsilon_{abc}\xi_c$$
which implies the integrability. Note that the Reeb distribution
is the vertical distribution of the Konishi bundle.

 \noindent Let
$g_{\OO}$ be the metric on $\OO$, and define the following
pseudoriemannian metric (with signature $(4k+3,4l)$) on $\SS$:
$$g_{\SS}:=\pi^*g_{\OO}+\sum_i\eta_i\otimes\eta_i.$$
Since $\xi_i$ are unit Killing vector fields w.r.t. $g_{\SS}$, the
tensors $$\Theta_i:=\nabla\xi_i$$ are skewsymmetric, and they
satisfy
\begin{itemize}
\item $\Phi_i^2=-{\rm Id}+\eta_i\otimes\xi_i$,%
\item $g_{\SS}\circ(\Phi_i\otimes\Phi_i)=g_{\SS}-\eta_i\otimes\eta_i$,%
\item $g_{\SS}\circ(\Phi_i\otimes {\rm Id})=d\eta_i$.
\end{itemize}
In particular, $(\eta_i,\xi_i,\Phi_i)$ are \emph{contact metric
structures} w.r.t. $g_{\SS}$. Consider the \textbf{metric cone} on
$\SS$ (sometimes called the \textbf{Swann's bundle} \cite{Sw91} of
$\OO$ and denoted by $\UU(\OO)$), i.e.
$$(\CC(\SS),g_{\CC}):=(\SS\times\R^+,t^2g_{\SS}+dt^2).$$
The cone $\CC(\SS)$ is endowed with three almost complex
structures
$$J_i:\left\{\begin{array}{lcl}J_i(X)&=&\Phi_i(X)-\eta_i(X)\Psi\\
                                     &&\\
                               J_i(\Psi)&=&\xi_i\end{array}\right.,\quad X\in\Gamma(T\SS),\quad \Psi=t\partial_t.$$
Indeed, $(J_1,J_2,J_3,g_{\CC})$ is a \textbf{pseudo-hyperk\"ahler
structure} (with signature $(4k+4,4l)$) on $\CC$. In other words,
$(\{\eta_i,\xi_i,\Phi_i\}_{i=1}^3,g_{\SS})$ is a
\textbf{pseudo-3-Sasakian structure} on $\SS$. In particular, the
metric $g_{\SS}$ is Einstein.
\subsection{The twistor space}
Pick any element $\tau\in S^2\subset\mathfrak{so}(3)$. The circle
$$S^1_{\tau}\subset{\rm SO}(3)$$ generated by $\tau$ acts freely and
isometrically on $\SS$, so we can define the \textbf{twistor
space} \cite{Sal82} $(\ZZ_{\tau},g_{\ZZ_{\tau}})$, where
$\ZZ_{\tau}:=S^1_{\tau}\setminus\SS$ and $g_{\ZZ_{\tau}}$ is the
metric (with signature $(4k+2,4l)$) obtained
from $g_{\SS}$ via the quotient.\\
If $\tau=(a_1,a_2,a_3)$, define
$$\eta_{\tau}:=a_1\eta_1+a_2\eta_2+a_3\eta_3,$$
$$\Phi_{\tau}:=a_1\Phi_1+a_2\Phi_2+a_3\Phi_3,$$
$$\DD_{\tau}:=\ker \eta_{\tau}.$$
The pair $(\DD_{\tau},\Phi_{\tau})$ defines an
$S^1_{\tau}$-invariant \emph{CR-structure} on $\SS$, hence an
almost complex structure $I_{\tau}$ on $\ZZ_{\tau}$. In fact, it
turns out that $(\ZZ_{\tau},I_{\tau},g_{\ZZ_{\tau}})$ is
pseudo-K\"ahler.\\
Finally, if $\tau=\tau_1\times\tau_2$ we can define the following
$\C$-valued 1-form on $\SS$, namely
$$\Upsilon_{\tau}^{{\tau}_1,{\tau}_2}:=\eta_{{\tau}_1}+i\eta_{{\tau}_2}.$$
Furthermore, the complex line bundle $L_{\tau}$ generated by
$\Upsilon_{\tau}^{{\tau}_1,{\tau}_2}$ doesn't depend on
$\tau_1,\tau_2$ and it is $S^1_{\tau}$-invariant, so it pushes
down to a bundle
$$\LL_{\tau}\subset T^{(1,0)}\ZZ_{\tau}$$
which defines a \textbf{pseudo-holomorphic contact structure} on
$\ZZ_{\tau}$, whose isomorphism class doesn't depend on the choice
of $\tau\in S^2$. In particular, the biholomorphism class of
$\ZZ_{\tau}$ doesn't depend on $\tau\in S^2$.\\
\\
Putting all these structures together, we obtain the following
diagram:
$$
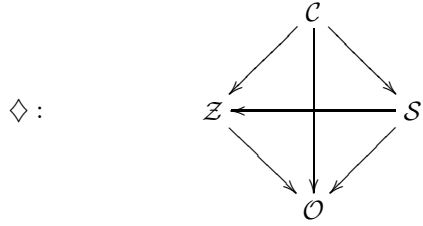
\xymatrix{&&&\CC\ar[dr]\ar[dl]\ar[dd]&\\
                  \diamondsuit:&&\ZZ\ar[dr]&&\ar[ll]\SS\ar[dl]\\
                   &&&\OO&}$$
\section{The reduced diamond diagram}
\noindent Let $\OO$ be a QK manifold with a group of isometries
$T$ preserving  the bundle $\mathfrak{H}$. Relative to this
action, we can associate the notion of a \textbf{moment map} (see
\cite{Ga87CMP} and \cite{GaLa88}), which is a \emph{section}
$$\mu:\OO\rightarrow \gt^*\otimes\mathfrak{H},$$
$$x\mapsto \{v\mapsto \mu_v(x)\},$$
where $\gt$ is the Lie algebra of $T$. The moment map is has  the
following properties:\begin{itemize} \item it is $T$-invariant,
i.e.
$$\mu_{{\rm Ad}_tv}(tx)=\mu_{v}(x)$$ \item it satisfies the
equation
$$\nabla\mu_v=i_{v_{\OO}}\Theta,$$
where $v_{\OO}$ is the vector field induced on $\OO$ by $v\in\gt$.
\end{itemize}
The moment map \emph{always exists and it is unique} as it is
completely determined by the second of the above conditions. In
fact, we have (\cite{Ga87CMP}):
$$\mu_v=\frac{\lambda'}{s}\sigma({\LL}_{v_{\OO}}-\nabla_{v_{\OO}}),$$
where $\sigma$ is the bundle isomorphism
$$\sigma:{\rm SkewEnd}(\mathfrak{H})\rightarrow \mathfrak{H}$$
such that
$$\sigma^{-1}(J)=\frac{1}{2}{\rm ad}_J=J\wedge.$$
It can be proven that, under appropriate assumptions on the action
and on the regularity of $\mu$, the quotient
$T\Qquot\OO:=T\setminus\mu^{-1}(0)$ is a QK orbifold (see
\cite{GaLa88}). The whole process
$$\xy\xymatrix{\OO\ar@{.>}[r]^{T\ \ \ \ }="A" &T\Qquot\OO}\endxy$$
is then called QK-\emph{reduction} of $\OO$.
\subsection{Diamond reductions}
Let $\phi:\OO\rightarrow\OO$ be a diffeomorphism preserving the QK
structure $(\gH,g_{\OO})$, i.e.,
$$\phi^*g=g,\quad\phi^*\gH=\gH.$$
Given a local frame $\sigma=\{J_1,J_2,J_3\}\in \Gamma(\SS)$, we
have $$\phi_* J_i:=(\phi^{*})^{-1} J_i:=\phi_*\circ
J_i\circ\phi_*^{-1}=\beta_i^jJ_j,$$ where
$\beta_{\sigma}:=\{\beta_i^j\}$ is a locally defined map with
values in ${\rm SO}(3)$, so $\phi$ induces a bundle automorphism
which we denote by $\phi_S$. There exists a map
$$\beta_{\phi}:\SS\rightarrow{\rm SO}(3)$$
such that $\beta\circ\sigma=\beta_{\sigma}$ for any local frame
$\sigma$ or, in other words,
$$\phi_{\SS}(f)={(\beta_{\phi}(f))}\cdot f.$$
Furthermore, if $\alpha_{\sigma}=\{\alpha_i^j\}$ is defined by
$\nabla J_i=\alpha_i^j\otimes J_j$, it follows that
$$\alpha_{\sigma}=\phi^*(\alpha_{{(\beta_{\phi}(\sigma))}\cdot\sigma})=\phi^*(\alpha_{\phi_{\SS}\circ\sigma}).$$
Hence, $$\phi_{\SS}^*\eta=\eta$$ and
$$\phi_{\SS}^*g_{\SS}=g_{\SS}.$$
Since $\phi_{\SS}$ preserves the $3$-Sasakian structure, it is
easy to see that it induces an automorphism of the twistor space.
Finally, the automorphism of the cone is obtained by trivially
lifting $\phi_{\SS}$.\\
Let $\pi:P\rightarrow\OO$ be any one of the three bundles in the
diamond diagram. If $\mu$ is the QK-moment map on $\OO$ and $X$ is
a Killing vector field on $P$, define a map
$$\widetilde{\mu}_X:=\mu_{\pi_*X}\circ\pi$$
which is a section of the pull-back bundle $\pi^*\gH$.\\

\noindent From the axioms of the QK moment map, it follows that
$$\nabla\widetilde{\mu}_X=i_X\pi^*\Theta$$
\begin{lem}For $P=\SS$ or $\CC$, the bundle $\pi^*\gH$ is trivial. Moreover,\begin{itemize}
\item if $P=\SS$, then $\widetilde{\mu}=\eta$;%
\item if $P=\CC$, then $\widetilde{\mu}$ is an hyperk\"ahler
moment map.
\end{itemize}
\end{lem}

\noindent The map $\eta$ is called \textbf{3-Sasakian moment map},
and the map induced on $\ZZ$ (for any fixed Killing field) is  a
section of the
bundle $\R\oplus\LL$ and is called \textbf{twistor moment map}.\\
We conclude the following: If $T$ is a group of QK isometries
acting on $\OO$ then it induces an action on the whole diamond
diagram  preserving all the relevant structures, and the QK moment
map can be extended to a moment map of the whole $\diamondsuit$.
Consequently it makes more sense to talk about the \textbf{diamond
reduction}:

$$\xy\xymatrix{\diamondsuit\ar@{.>}[r]^{T\ \ \ \ }="A" &T\Qquot\diamondsuit}\endxy.$$

\noindent It can be proven that, under appropriate assumptions,
the reduced diamond diagram $T\Qquot\diamondsuit$ is the diamond
diagram of the QK orbifold $T\Qquot\OO$ \cite{Sw91, BGM94a}.
\section{Wolf spaces}
\noindent The classification of (pseudo-)riemannian symmetric
spaces reduces to the study of involutive Lie algebras, i.e., Lie
algebras endowed with a canonical splitting (of vector subspaces)
$$\gg=\gh\oplus\gm,$$
where $\gh$ is a subalgebra and
$$[\gh,\gm]\subset\gm\quad [\gm,\gm]\subset\gh$$
In particular, any involutive Lie algebra defines corresponding
homogeneous space $M=G/H$, where $G,H$ are connected Lie groups
(generated by $\gg,\gh$ respectively), and $G$ is simply
connected.\\
The canonical splitting defines an $\ell_G$-invariant connection
on the principal bundle $\pi:G\rightarrow G/H$, given by
$$T_gG=l_{g*}\gh\oplus l_{g*}\gm$$
Via the horizontal lifting, for any open set $V\subset M$ we
construct a linear isomorphisms between sections and equivariant
maps
$$\Gamma_V(TM)\rightarrow C_H^{\infty}(\pi^{-1}(V),\gm)$$
$$X\rightarrow X^*.$$
In particular, $TM=G\times_H\gm$. We have similar isomorphisms for
other types of tensors: for example, the set of all the metrics
which give $M$ the structure of a pseudoriemannian symmetric space
is the set $\Sigma$ of all the ${\rm Ad}_H$-pseudoscalar products
on $\gm$.
\begin{teo}Let $\nabla$ be the linear connection on $TM$ associated to the canonical splitting of $\gg$. Then
$\nabla$ is $l_G$-invariant and is the Levi-Civita connection of
$TM$ with respect to all the elements of $\Sigma$. Furthermore,
its curvature two-form $R$ is associated to a map
$$R^*:G\rightarrow \bigwedge^2\gm^*$$
given by
$$R^*(v,w)=-[v,w]$$ \end{teo}
Thus, by the Ambrose-Singer theorem, the holonomy algebra is given
by ${\rm ad}_{[\gm,\gm]|\gm}$. The irreducibility of this algebra
implies that $\gg$ is semisimple. As a scalar
product we can thus take the restriction of the Killing form of $\gg$ to $\gm$.\\
Motivated by all these facts we give the following
\begin{defn} A [generalized] \emph{Wolf space} is a semisimple, involutive
Lie algebra $\gg=\gh\oplus\gm$ such that
$$[\gm,\gm]=\gh=\gh'\oplus\gs,$$
where the last splitting is the direct sum of two commuting Lie
subalgebras, such that $\gs\simeq\mathfrak{sp}(1)$ and $\gh'$ is
isomorphic to a Lie subalgebra of $\mathfrak{sp}(l,k)$ with
$4k+4l={\rm dim}\ \gm$.
\end{defn}
\noindent It follows that a Wolf space is
(pseudo-)quaternion-K\"ahler, and that the adjoint representation
of $H=H'S$ (where $\gs={\rm Lie}(S)$, $\gh'={\rm Lie}(H')$) on
$\gm$ and the usual one are isomorphic. Moreover, the action of
$G$ on $M$ is almost effective and the Killing algebra of $M$ is
$\gg$.  The bundle
$$\mathfrak{H}:=G\times_H\gs$$
defines the almost quaternionic hermitian structure, and we have
the linear isomorphisms
$$\Gamma_V(\mathfrak{H})\rightarrow C^{\infty}_H(\pi^{-1}(V), \gs)$$
$$J\mapsto J^*$$
such that
$$(JX)^*=[J^*,X^*],$$
as long as $X$ is left invariant. Moreover, the linear connection
induced on $\mathfrak{H}$ is the restriction  of the
Levi-Civita connection on tensors to $\mathfrak{H}$.\\
The 3-Sasakian bundle is given by
$$\SS=G\times_H{\rm SO}(3)=G/H'\times_{{\rm SO}(3)}{\rm SO}(3)=G/H'$$
with the trivial projection. The $1$-form defining the 3-sasakian
structure is the connection $1$-form associated to the Levi-Civita
connection of $\mathfrak{H}$, i.e.,
$$\eta(l_{g*}v)=v_{\gs}$$
for each $v\in \gg$.\\
Now, pick $i\in S^2\subset\gs$. This induces a vector field and an
$U$-action on $\SS$, where $U\simeq {\rm U(1)}$ is the centralizer
of $i$ in $S$ and the action is nothing but the right
multiplication. Therefore,
$$\ZZ=G/H'U$$
In conclusion, the diamond diagram is given by
$$\xymatrix{&{\CC}(G/H')\ar[dr]\ar[dl]\ar[dd]&\\
                     G/H'U\ar[dr]&&\ar[ll]G/H'\ar[dl]\\
                    &G/H&}$$

\subsection{The Alekseevsky-Cort{\'e}s's list}

\noindent Recently, the generalized Wolf spaces have been
classified (see \cite{AlCor05}):
$$\frac{{\rm SU}(p+2,q)}{{\rm S}({\rm U}(2)\times {\rm U}(p,q))}\quad \frac{{\rm SL}(n+1,\H)}{{\rm S}({\rm GL}(1,\H)\times{\rm GL}(n,\H))}$$
$$\frac{{\rm SO}_o(p+4,q)}{{\rm SO}(4)\times {\rm SO}_o(p,q)}\quad \frac{{\rm SO}^*(2l+4)}{{\rm SO}^*(4)\times {\rm SO}^*(2l)}$$
$$\frac{{\rm Sp}(p+1,q)}{{\rm Sp}(1)\times {\rm Sp}(p,q)}$$
$$\frac{E_{6(-78)}}{{\rm SU}(2){\rm SU}(6)}\quad \frac{E_{6(2)}}{{\rm SU}(2){\rm SU}(6)}\quad \frac{E_{6(2)}}{{\rm SU}(2){\rm SU}(2,4)}\quad \frac{E_{6(-14)}}{{\rm SU}(2){\rm SU}(2,4)}$$
$$\frac{E_{6(6)}}{{\rm Sp}(1){\rm SL}(3,\H)}\quad \frac{E_{6(-26)}}{{\rm Sp}(1){\rm SL}(3,\H)}$$
$$\frac{{\rm E}_{7(-133)}}{{\rm SU}(2){\rm Spin}(12)}\quad \frac{{\rm E}_{7(-5)}}{{\rm SU}(2){\rm Spin}(12)}\quad \frac{{\rm E}_{7(-5)}}{{\rm SU}(2){\rm Spin}_o(4,8)}\quad \frac{{\rm E}_{7(7)}}{{\rm SU}(2){\rm SO}^*(12)}\quad \frac{{\rm E}_{7(-25)}}{{\rm SU}(2){\rm SO}^*(12)}$$
$$\frac{{\rm E}_{8(-248)}}{{\rm SU}(2)E_{7(133)}}\quad \frac{{\rm E}_{8(-24)}}{{\rm SU}(2)E_{7(133)}}\quad \frac{{\rm E}_{8(-24)}}{{\rm SU}(2)E_{7(-5)}}\quad \frac{{\rm E}_{8(8)}}{{\rm SU}(2)E_{7(-5)}}$$
$$\frac{{\rm F}_{4(-52)}}{{\rm Sp}(1){\rm Sp}(3)}\quad \frac{{\rm F}_{4(4)}}{{\rm Sp}(1){\rm Sp}(3)}\quad \frac{{\rm F}_{4(4)}}{{\rm Sp}(1){\rm Sp}(1,2)}\quad \frac{{\rm F}_{4(-20)}}{{\rm Sp}(1){\rm Sp}(1,2)}$$
$$\frac{{\rm G}_{2(-14)}}{{\rm SO}(4)}\quad \frac{{\rm G}_{2(2)}}{{\rm SO}(4)}$$
In the compact Riemannian case, we get the classical \emph{classic
Wolf spaces}, which arise from the classification of all QK
homogeneous spaces with positive scalar curvature (see
\cite{Wol65}), namely the spaces
$$\frac{{\rm Sp}(n+1)}{{\rm Sp}(n)\times {\rm Sp}(1)}\quad\frac{{\rm SU}(n)}{{\rm S}({\rm U}(n-2)\times{\rm U}(1)){\rm Sp}(1)}\quad\frac{{\rm SO}(n)}{{\rm SO}(n-4)\times {\rm SO}(4)}$$
$$\frac{G_2}{{\rm
SO}(4)}\quad \frac{F_4}{{\rm Sp}(3){\rm Sp}(1)}\quad
\frac{E_6}{{\rm SU }(6){\rm Sp}(1)}\quad \frac{E_7}{{\rm
Spin}(12){\rm Sp}(1)}\quad \frac{E_8}{E_7{\rm Sp(1)}}$$ It is
straightforward to verify that the set of generalized Wolf spaces
is closed with respect to the duality of involutive Lie algebras,
i.e. the map
$$\gg=\gh\oplus\gm\ \longmapsto\ \gh\oplus i\gm\subset\gg\otimes\C$$
\section{Quaternionic reductions of Wolf spaces}

\subsection{Formula for the moment map for Wolf spaces}

\noindent Suppose we have now a Wolf space of the form $G/H'S$,
and that the action is the left multiplication by $T$, where $T$
is some \emph{vitual} Lie subgroup of $G$ (i.e., a subgroup of
$G$, not necessarily closed). Then, it can be easily seen that this action preserves all the structures in the diamond diagram.\\

\noindent The aim of this work is to study (and eventually give
some classification) of all the possible reductions we could get
in this way, when $T$ varies among all the virtual Lie subgroups
of $G$.\\

\noindent Actually, in order to produce Hausdorff spaces, we'll
require a \emph{proper action} and in particular that $T$ is
indeed a Lie subgroup. For the moment we don't care and focus our
attention on the Lie algebra which \emph{generates the action},
i.e.
$$\gt:={\rm Lie }(T)\subset \gg$$

\noindent The adventage of working with (generalized) Wolf spaces
is quite clear: The moment map has a particularly simple form. Up
to some scalar factor, it is given by the formula
$$\mu_v(g)=({\rm Ad}_g^{-1}v)_{\gs}$$
(seen as an equivariant map), for any $g\in G$ and $v\in \gt$.\\
(Note that the quaternionic $\mathfrak{H}$-valued two-form is
given by
$$\Theta(v,w)=-\lambda[v,w]_{\gs}$$
for some $\lambda>0$.)\\
As a matter of fact, this map \emph{at the same time} represents
the 3-sasakian and the quaternion-K\"ahler moment map, so we'll
call it just \emph{moment map}.
\subsection{The reduced diamond diagram}

\noindent We can define four \emph{zero loci}:
$$Z_G(T):=\{g\in G:\mu_v(g)=0 \mbox{ for each }v\in \gt\},$$
$$Z_{\SS}(T):=\{gH'\in \SS:\mu_v(g)=0 \mbox{ for each }v\in \gt\}=Z_G/H',$$
$$Z_{\ZZ}(T):=\{gH'U\in \ZZ:\mu_v(g)=0 \mbox{ for each }v\in \gt\}=Z_G/H'U,$$
$$Z_{\OO}(T):=\{gH\in \OO:\mu_v(g)=0 \mbox{ for each }v\in \gt\}=Z_G/H,$$
and the reductions
$$\RR_G(T):=T\setminus Z_G(T),\quad\RR_{\SS}(T):=T\setminus Z_{\SS}(T),\quad\RR_{\ZZ}(T):=T\setminus Z_{\ZZ}(T),\quad\RR_{\OO}(T):=T\setminus Z_{\OO}(T),$$
giving a diagram
$$
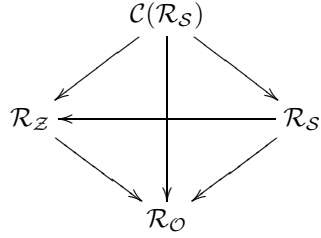
\xymatrix{&{\CC}({\RR_{\SS}})\ar[dr]\ar[dl]\ar[dd]&\\
                    {\RR_{\ZZ}}\ar[dr]&&\ar[ll]{\RR_{\SS}}\ar[dl].\\
                    &{\RR_{\OO}}&}$$
Our formula of the moment map implies an \emph{extra symmetry}
with respect to the group $G$ action, namely
$$Z_G(T^g)=g\cdot Z_G(T).$$
Hence, we can restrict our classification only to \emph{conjugacy
classes} of virtual Lie subgroups of $G$. Furthermore, it follows
that the group $C_G(T)/T$ preserves the reduced diamond diagram.

\begin{quest}
Is this group \emph{the stabilizer} of the diagram?
\end{quest}

\noindent There are a lot of problems in the reduction process.
The reduced spaces may not have any orbifold structure  or even
fail to be Hausdorff. This, for example, can be due to some
irregularity of the moment map or non-finiteness of some isotropy
subgroups of $T$, the noncompactness of $T$, etc. Nevertheless, we
shall examine examples (perhaps exceptional), where the moment map
is not regular, some isotropies are not finite and one still can
give to these spaces an orbifold structure in such a way that the
above diagram is actually a diamond diagram.

\section{1-dimensional toric reductions}
\subsection{Singular and irregular points}

\noindent Let us analyze the action of $T$ on the zero loci. Note
that tt is automatically free on $Z_G$ but on the other spaces we
could get nontrivial isotropies. For instance, $tgK=gK$ if and
only if $t\in T\cap K^g$ and in our cases we get three different
kinds of isotropies (respectively, \emph{3-Sasakian},
\emph{twistor} and \emph{quaternionic}), namely
$$T\cap (H')^g\subset T\cap (H'U)^g\subset T\cap (H)^g.$$
For each type of action we distinguish two different kinds of
points with nontrivial isotropy, namely the \emph{singular points}
(with discrete isotropy) and \emph{irregular points} (all the
other), so we have 3-Sasakian irregular and singular points...and
so on.\\
\begin{lem}If $T$ is 1-dimensional, then $Z_G$ contains only irregular points of 3-Sasakian type.\end{lem}
\emph{Proof}. Since $T$ is 1-dimensional, it doesn't have
nondiscrete, proper subgroups. Hence, the three irregularity
conditions can be written in the following way:
$${\rm Ad}_{g}^{-1}\gt\subset\gh',$$
$${\rm Ad}_{g}^{-1}\gt\subset\gh'\oplus\gu,$$
$${\rm Ad}_{g}^{-1}\gt\subset\gh.$$
Comparing these conditions with the moment map equation
$${\rm Ad}_{g}^{-1}\gt\subset\gh'\oplus\gm$$
proves the lemma.\qed\\
For this reason, by \emph{irregular points} we shall mean just
irregular points of 3-Sasakian type.
\subsubsection{The regularity condition}
Now, let's write down the differential of the moment map: it can
be thought as a map
$$d\mu_v:G\times\gg\rightarrow\gs,$$
$$(g,w)\rightarrow d_g\mu_v(w).$$
An easy calculation shows that
$$d_g\mu_v(w)=-\mu_{[{\rm Ad}_gw,v]}(g)=-({\rm Ad}_g^{-1}[{\rm Ad}_gw,v])_{\gs}=[{\rm Ad}_g^{-1}v,w]_{\gs}$$
\begin{lem}If $T$ is 1-dimensional then on $Z_G$ the critical set of the moment map coincides with its irregular set
and with the set of (3-Sasakian) irregular points.\end{lem}
\emph{Proof:} Let $<,>$ be the Killing form on $\gg$, $w\in{\rm
Ker}\ d_g\mu_v$ for some $g\in Z_G$. Then, for each $\xi\in\gg$ we
have
$$0=<d_g\mu_v(w),\xi>=<[{\rm Ad}_g^{-1}v,w]_{\gs},\xi>=-<w,[{\rm Ad}_g^{-1}v,\xi_{\gs}]>.$$
Hence, ${\rm Ker}\ d_g\mu_v$ is the orthogonal complement of
$[{\rm Ad}_g^{-1}v,{\gs}]$ which has codimension $<3$ if and only
if
$({\rm Ad}_g^{-1}v)_{\gm}=0$.\qed\\

\noindent In particular, this lemma says that, if the $3$-Sasakian
action is locally free, the normal bundle of $Z_G$ in $G$ is
$$\nu(Z_G):=\{(g,\xi)\in S\times\gg:g\in Z_G,\xi\in [{\rm Ad}_g^{-1}v,{\gs}]\}\simeq Z_G\times \gs.$$
Thus, $Z_G$ is parallelizable. Moreover, the orthogonal projection
$$\pi:TG_{|Z_G}\rightarrow\nu(Z_G)$$
can be thought as a map
$$\pi:\gg\rightarrow\gs$$
satisfying
$$\pi(w)=\frac{\lambda}{\|({\rm Ad}_g^{-1}v)_{\gm}\|}d_g\mu_v(w)$$
provided that the denominator doesn't vanish.\\
Again, $\pi$ can be pushed down to a bundle map
$$\pi':T\OO_{|Z_{\OO}}\rightarrow\nu(Z_{\OO}).$$
\begin{lem}The \emph{second fundamental form} associated to the embedding $Z_{\OO}\subset\OO$ is given by the equivariant map
$$\alpha:\alpha_g(w_1,w_2)=\frac{\lambda}{\|({\rm Ad}_g^{-1}v)_{\gm}\|}[[w_1,{\rm Ad}_{g}^{-1}v],w_2]_{\gs}.$$
\end{lem} \emph{Proof}: Let $X,Y:Z_G\rightarrow\gm$ be vector fields of $Z_{\OO}$.  The covariant derivative
$\nabla_XY$, thought as equivariant map, is given by
$$(\nabla_XY)_g=(d_gY)(l_{g*}X_g).$$
Thus, the second fundamental form is given by
$$\pi((d_gY)(l_{g*}X_g))=\frac{\lambda}{\|({\rm Ad}_g^{-1}v)_{\gm}\|}[{\rm Ad}_g^{-1}v,(d_gY)(l_{g*}X_g)]_{\gs}.$$
As
$$[{\rm Ad}_g^{-1}v,Y_g]_{\gs}=0\ \forall\ g\in Z_G$$
by taking $l_{g*}X_g\in T_gZ_G$, we get
$$0=-[[X_g,{\rm Ad}_g^{-1}v],Y_g]_{\gs}+[{\rm Ad}_g^{-1}v,(d_gY)(l_{g*}X_g)]_{\gs}.$$
\qed\\
We now are able to write the sectional curvature of $Z_{\OO}$: it
is given by
\begin{eqnarray}\mathfrak{K}_g(w_1,w_2)&=&\|[w_1,w_2]\|+<\alpha_g(w_1,w_1),\alpha_g(w_2,w_2)>-\|\alpha_g(w_1,w_2)\|,\nonumber\end{eqnarray}
where the vectors $w_1,w_2$ are orthogonal.

\subsection{The quaternion-K\"ahler structure of $\RR_{\OO}$}

\noindent Suppose that the group $T$ acts on $Z_G$ with finite
isotropy. Then, the bundle
$$T\setminus\nu(Z_G)/H=\RR_G\times_H\gs$$
defines the quaternionic structure on $\RR_{\OO}$: in fact, it has
fiber $\simeq\gs$ and can be thought as a subbundle of ${\rm
End}(T\RR_{\OO})$. Furthermore, the pull-back of the metric of
$\OO$ on $Z_{\OO}$ can be pushed to a metric on $\RR_{\OO}$.

\subsection{The energy function}

\noindent Related to the moment map we have the map
$$E:G\rightarrow [0,+\infty),$$
$$\ \ g\mapsto-\|\mu_v(g)\|,$$
whose gradient is given by
$$({\rm grad}\ E)_g=-[({\rm Ad}_g^{-1}v)_{\gs},({\rm Ad}_g^{-1}v)_{\gm}].$$
Its critical set is given by the union of $Z_G$ and all the
quaternionic irregular points in $G$. It turns out that this
vector field is \emph{complete} (see \cite{Bat04}), so it defines
a global flow on $G$. Furthermore, if the adjoint orbit of $v$
doesn't intersect $\gh$ (note that, if $H$ is compact, at most one
adjoint orbit satisfies this condition, while it doesn't exist if
$G$ is compact),
$Z_G$ is the critical set of $E$.\\
Unfortunately, we do not know if the gradient flow converges to
$Z_G$ when $G$ is not compact:
\begin{quest}
Is the zero locus $Z_G$ (resp. $Z_{\SS}$, $Z_{\UU}$, $Z_{\OO}$) a
deformation retract of $G$ (resp. $\SS$, $\UU$, $\OO$)?
\end{quest}

\section{${\rm SO}(7)/{\rm SO}(3)\times {\rm SO}(4)$ and ${\rm SO}_o(3,4)/{\rm SO}(3)\times {\rm SO}(4)$}

\noindent Let us consider the spaces ${\rm SO}(7)/{\rm
SO}(3)\times {\rm SO}(4)$ and ${\rm SO}_o(3,4)/{\rm SO}(3)\times
{\rm SO}(4)$. The first one is the usual Grassmann manifolds of
4-planes in $\R^7$, or in other words the space of all the
orthogonal splittings
$$W^4\oplus (W^4)^{\bot}=\R^7$$
with respect to the usual euclidean scalar product.\\
The latter space is (the connected component of some fixed point
$o$ of) the space of splittings as before, but with respect the
standard \emph{pseudo-euclidean scalar product with signature}
$(3,4)$:
$$ds^2=dx_1^2+dx_2^2+dx_3^2-dx_4^2-dx_5^2-dx_6^2-dx_7^2,$$
and such that $W^4$ and $(W^4)^{\bot}$ are \emph{space-like} and
\emph{time-like} subspaces, respectively (i.e., such that the
restriction of the scalar product is respectively negative and
positive defined).\\
We shall study all the possible homogeneous 1-dimensional
quaternionic K\"ahler reductions. In order to do this, it will be
useful to give a list of all the adjoint orbits of the group $G$.
Actually, since the first group is compact the classification is
elementary: all the adjoint orbits in ${\rm SO}(7)$ have a
representative of the form
$$\left(\begin{array}{ccccccc}0&0&0&0&0&0&0\\
                              0&0&a&0&0&0&0\\
                              0&-a&0&0&0&0&0\\
                              0&0&0&0&b&0&0\\
                              0&0&0&-b&0&0&0\\
                              0&0&0&0&0&0&c\\
                              0&0&0&0&0&-c&0\end{array}\right)$$
where $a,b,c$ are nonnegative real numbers.\\
The classification of adjoint orbits of ${\rm SO}(3,4)$ is much
more complicated and gives several families of adjoint orbits,
each of them depending on some set of parameters. In order to do
that, we'll quote an algorithm found by Burgoyne and Cushman (see
\cite{BurCush}), valid
for all classical groups.\\
\subsection{Octonions and split octonions}

\noindent Before to start the classification of adjoint orbits, I
want to consider for both of this spaces a particular reduction,
which involves some considerations about \emph{normed algebra}
structures on $\R^8$.
\subsubsection{The Cayley-Dickson process}

\noindent (See \cite{Harv90} for more details.) We first recall
what octonions and split octonions are, starting
from some general facts about normed algebras.\\
A normed algebra is a (real) algebra with unity $V$ endowed with
some non-degenerate quadratic form $\|\cdot\|$ such that
$\|vw\|=\|v\|\|w\|$. Via polarization, this quadratic form can be
identified with some pseudo-scalar product $<,>$ such that
$$<xz,yw>+<xw,yz>=2<x,y><z,w>.$$
It follows that $\|1_V\|=1$, so ${\rm Re}(V):=\R 1_V$ has an
orthogonal complement in $V$ which we denote with ${\rm Im}(V)$.
In a normed algebra we can define a \emph{conjugation} by using
this splitting in a natural way.\\
Now, given a  normed algebra $(V,\|\cdot\|)$ we can construct two
new algebras $V(\pm)$, whose underlying vector space is nothing
but $V\oplus V$ and which contain $V$ as a subalgebra. This
construction (called \emph{Cayley-Dickson process}) is motivated
by the following
\begin{lem}If $A$ is a normed subalgebra (with $1\in A$) of a normed algebra $B$
and $\epsilon \in A^{\bot}$ is such that $\|\epsilon\|=\pm 1$,
then $\epsilon A$ is orthogonal to $A$ and
$$(a+\epsilon b)(c+\epsilon d)=(ac\mp \overline{d}b)+\epsilon (da+b\overline{c})$$\end{lem}

\noindent Hence, the product on $V(\pm)$ must be
$$(a, b)(c, d)=(ac\mp \overline{d}b,da+b\overline{c})$$
while the new norm is given by
$$\|(a,b)\|:=\|a\|\pm\|b\|$$
The Cayley-Dickson process has the following properties:
\begin{itemize}
\item $V(\pm)$ is commutative if and only if $V=\R$; \item
$V(\pm)$ is associative if and only if $V$ is commutative and
associative; \item $V(\pm)$ is normed iff $V(\pm)$ is alternative
(i.e., the associator $[x,y,z]$ is alternating), iff $V$ is
associative.
\end{itemize}
In particular, if we start from $\R$ we get the following diagram:
$$
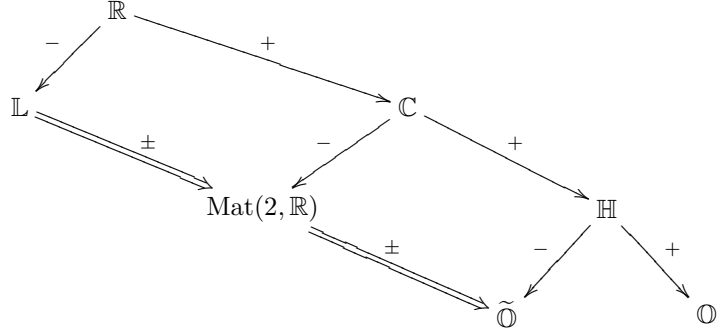
\xymatrix{&\R\ar@{->}[dl]_{-}\ar@{->}[drr]^+&&&&& \\
            {\mathbb L}\ar@{=>}[drr]^{\pm}&&&\C\ar@{->}[dl]_{-}\ar@{->}[drr]^+&&& \\
            &&{\rm Mat}(2,\R)\ar@{=>}[drr]^{\pm}&&&\H\ar@{->}[dl]_{-}\ar@{->}[dr]^+& \\
            &&&&\widetilde{\O}&&\O }$$
After the third step, this process doesn't produce other normed
algebras. The well-known theorem due by Hurwitz says that
these are the only possible normed algebras.\\
In particular, we want to focus our attention on the last two
algebras $\O$ (the \emph{octonions}), $\widetilde{\O}$ (the
\emph{split octonions}) and their groups of automorphisms.\\
Let's consider first the quaternions
$\H=<1,I_1,I_2,I_3>\simeq\R^4$ and define
$$1=e_0=(1,0),\quad e_4=(0,1)$$
$$e_k=(I_k,0),\quad e_{k+4}=(0,I_k)$$
The algebra structure of the octonions is completely described by
a diagram (called \textbf{the Fano plane}), that with respect to
the basis $e_1,e_2,e_3,e_4,e_5,e_6,e_7$ has the
form\xyoption{frame}
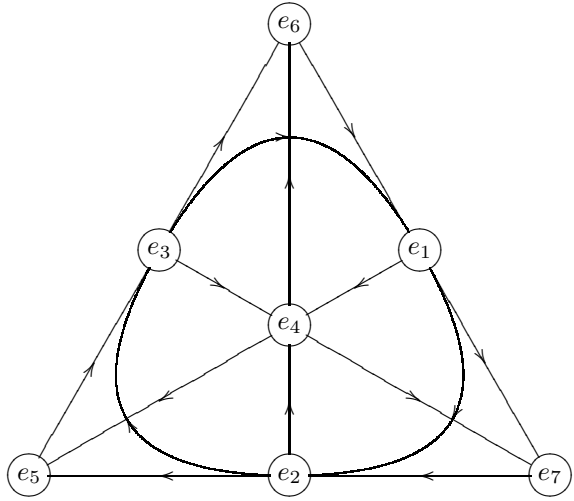
\xy <9.5cm,-2.5cm>;<11.5cm,-2.5cm>:(0,0)+*+<1.33mm,1.33mm>+/5.1pt/{e_4}="4"*\frm<8pt>{o};%
(0,1)+*+<1.33mm,1.33mm>+/4.1pt/{e_6},*\frm<8pt>{o}="6";%
(0,-.5)+*+<1.88mm,1.88mm>+/5.1pt/{e_2},*\frm<8pt>{o}="2";%
(\halfrootthree,-.5)+*=\frm<8pt>{o},*+<1.8mm,1mm>+/5.1pt/{e_7}="7";%
(-\halfrootthree,-.5)+*+<1.8mm,1mm>+/5.1pt/{e_5},*\frm<8pt>{o}="5";%
(.43301270,.25)+*\frm<8pt>{o},*+<1.33mm,1.33mm>+/4pt/{e_1}="1";%
(-.43301270,.25)+*+<1.33mm,1.33mm>+/4pt/{e_3},*\frm<8pt>{o}="3";%
"1";"2"**\crv{"7"}?*\dir{>};"2";"3"**\crv{"5"}?*\dir{>};"3";"1"**\crv{"6"}?*\dir{>};%
"1";"4"**\dir{-}?*\dir{>};"1";"7"**\dir{-}?*\dir{>};%
"2";"4"**\dir{-}?*\dir{>};"2";"5"**\dir{-}?*\dir{>};%
"3";"4"**\dir{-}?*\dir{>};"3";"6"**\dir{-}?*\dir{>};%
"4";"5"**\dir{-}?*\dir{>};"4";"6"**\dir{-}?*\dir{>};"4";"7"**\dir{-}?*\dir{>};%
"7";"2"**\dir{-}?*\dir{>};"5";"3"**\dir{-}?*\dir{>};"6";"1"**\dir{-}?*\dir{>};%
\endxy

\ \\In fact, if $e_i\neq e_j$ belong to the same edge of $e_k$,
then $e_ie_j=\pm e_k$ where the sign depends on the order of
$e_i,e_j,e_k$ along the edge. Moreover, $e_i^2=-1$ for each $i$.\\
Furthermore, with respect to this basis we obtain the
multiplication rule of the split octonions by changing the sign to
the products of all elements in $\{a_4,a_5,a_6,a_7\}$. In the case
of octonions, the norm arising from the Cayley-Dickson process is
the euclidean one, while in the case of split octonions it has
signature $(4,4)$.
\begin{lem}Let $x,y,z\in {\O}$ (or $\in\widetilde{\O}$)  be purely imaginary and orthogonal.
Then \begin{itemize} \item $xy=-yx$ is purely imaginary and
orthogonal to $x,y$; \item the double product $x(yz)$ is
alternating.
\end{itemize}\end{lem}
\begin{lem}Let $x,y,z\in {\O}$ (or $\in\widetilde{\O}$).
Then we have the following $\text{Moufang identities}$:
\begin{itemize} \item $(xyx)z=x(y(xz))$;
\item $z(xyx)=((zx)y)x$; \item $(xy)(zx)=x(yz)x$.
\end{itemize}\end{lem}

\noindent Let $G_2,G_{2(2)}\subset {\rm GL}(8,\R)$ be the
automorphism groups of the algebras ${\O},\widetilde{\O}$,
respectively. Since they preserve the identity and the conjugation
they must be subgroups of ${\rm
O}(7)$, ${\rm O}(3,4)$, respectively.\\
>From now on we will use  $V$ to denote either $\R^7$ or
$\R^{3,4}$, and $<,>_V$ will be the (pseudo-)scalar product
associated to $V$. Moreover, we identify $V$ with the imaginary
part of $\O(V)$, where
$$\O(V):=\left\{\begin{array}{lcl}\O&\mbox{ if }&V=\R^7\\
                                    &&\\
                                  \widetilde{\O}&\mbox{ if }&V=\R^{3,4}\end{array}\right.$$
and $G_2(V)$ will be the automorphisms group of $\O(V)$. In
particular, $G_2(V)$ preserves the so-called \emph{associative
3-form} on $V$
$$\phi_V(x,y,z):=<x,yz>,\qquad \phi_V\in\bigwedge^3V^*$$
and the \emph{coassociative 4-form}
$$\psi_V(x,y,z,w):=<x,y(zw)-w(zy)>,\qquad \psi_V\in\bigwedge^4V^*$$
$G_2(V)$ preserves $\phi\wedge\psi=d{\rm vol}(V)$ and, hence, it
is a subgroup of ${\rm SO}(V)$.\\
Moreover, it can be shown that $G_2(V)$ is the stabilizer of
$\phi_V$ in ${\rm GL}(7)$. The equation
$$g^*\phi=\phi$$
implies that $G_2(V)$ is a 14-dimensional Lie group.\\ Taking the
derivative of the last relation, we obtain the relations which
define the Lie algebra $\gg_{2}(V)\subset\mathfrak{so}(V)$:
$$\left\{\begin{array}{lcl}a_{12}+a_{47}-a_{56}&=&0\\
                           a_{13}-a_{46}-a_{57}&=&0\\
                           a_{14}-a_{27}-\epsilon_V a_{36}&=&0\\
                           a_{15}+a_{26}-\epsilon_V a_{37}&=&0\\
                           a_{16}-a_{25}+\epsilon_V a_{34}&=&0\\
                           a_{17}+a_{24}+\epsilon_V a_{35}&=&0\\
                           a_{23}+a_{45}-a_{67}&=&0\end{array}\right.\qquad \left(\{a_{ij}\}\in \mathfrak{so}(V)\right)$$
where
$$\epsilon_V:=\left\{\begin{array}{lcl}+1&\mbox{ if }&V=\R^7,\\
                                    &&\\
                                       -1&\mbox{ if }&V=\R^{3,4}.\end{array}\right.$$
Furthermore, using the Moufang identities it is easy to prove
\begin{teo}A matrix $A=(a_1|a_2|a_3|a_4|a_5|a_6|a_7)\in {\rm GL}(7)$ belongs to $G_2(V)$ if and only if
$$\left\{\begin{array}{ccc}a_4a_5&=&\epsilon_V a_1\\
                           a_4a_6&=&\epsilon_V a_2\\
                           a_4a_7&=&\epsilon_V a_3\\
                           a_4a_5+a_6a_7&=&0\end{array}\right.$$
and $(a_4,a_5,a_6,a_7)$ is an orthogonal 4-frame with
$\|a_4\|=\|a_5\|=\|a_6\|=\|a_7\|=\epsilon_V$.
\end{teo}
\subsection{The canonical 1-dimensional reduction}

\noindent We are going to link the (split) octonions to a
particular
action on ${\rm SO}(V)$.\\
First of all, note that the scalar product on $V$ defines a linear
map
$$F:\mathfrak{so}(V)\rightarrow \bigwedge^2V^*$$
$$\ \ A\mapsto <A(\cdot),\cdot>$$
which is an isomorphism. In particular, every two-form $\omega$ on
$V$ defines a left action on ${\rm SO}(V)$ (and a QK action on
${\rm SO}(V)/{\rm SO}(3)\times{\rm SO}(4)$). In this case, the
moment map can be expressed in terms of the 2-form $\omega$:
$$\mu_{F^{-1}(\omega)}(g)=\left(\begin{array}{c}\omega(f_1,f_2)+\omega(f_3,f_4)\\
                  \omega(f_1,f_3)-\omega(f_2,f_4)\\
                  \omega(f_1,f_4)+\omega(f_2,f_3)\end{array}\right),$$
where $f_1,f_2,f_3,f_4$ are the last four columns of $g\in{\rm
SO}(V)$. In particular, the quaternionic zero locus can be
described as the set of all the euclidean 4-planes in $V$ such
that the restriction of $\omega$ on them satisfies
$*\omega=-\omega$.\\Now, fix $x\in V$ and set $\omega=i_x\phi_V$,
where $\phi_V$ is the associative 3-form on $V$.\\In this case,
the moment map equations are
$$<f_1f_2+f_3f_4,x>=0,$$
$$<f_1f_3-f_2f_4,x>=0,$$
$$<f_1f_4+f_2f_3,x>=0.$$
Hence we get
\begin{teo}The zero locus $Z_G^x$ of the action generated by $A_x:=F^{-1}(i_x\phi)$ contains the set $T_x\cdot G_2(V)\cdot{\rm SO}(3)$, where
$T_x$ is the group generated by $A_x$.\end{teo} \emph{Proof:}
$Z_G^x$ contains the set $\Sigma$ defined by equations
$$f_1f_2+f_3f_4=0,$$
$$f_1f_3-f_2f_4=0,$$
$$f_1f_4+f_2f_3=0.$$
The Moufang identities imply that each of these equations is
sufficient to describe $\Sigma$, and the statement follows from
the previous theorem.\qed\\We now search for the elements $x\in V$
such that the equality holds, i.e.
$$Z_G^x=T_x\cdot G_2(V)\cdot{\rm SO}(3)$$
>From the relation $A_{\lambda x}=\lambda A_x$ it follows that we
can restrict $x$ to
the sphere $S(V)\subset V$.\\
Moreover, for any $g\in G_2(V)$ we get
$$A_{g(x)}={\rm Ad}^{-1}_g(A_x).$$
Hence,
\begin{itemize}
\item if $V=\R^7$ we have just one case, since $G_2$ acts
transitively on $S^6$;%
\item if $V=\R^{3,4}$, $S(V)=S^+\cup S^-\cup K-\{0\}$, and each
component is an orbit of $G_{2(2)}$, so we have \emph{three}
cases, i.e. $x\in S^+$ (the \textbf{time-like case}), $x\in S^-$
(the \textbf{space-like case}) or $x\in K-\{0\}$ (the
\textbf{light-like case}).
\end{itemize}
Furthermore, $Z^x_G$ is acted on by the group $$U_x:=T_x\ltimes
H_x,$$ where $H_x$ is the isotropy subgroup of $x$ in $G_2(V)$.
\begin{teo}Suppose either $x\in S^6$ (in the compact case) or $x\in S^+$ (in the noncompact case).
Then, $$Z_G^x=T_x\cdot G_2(V)\cdot{\rm SO}(3)$$ and
$$\RR_G^x=\Z_3\setminus G_2(V)\cdot{\rm SO}(3).$$
\end{teo}
\noindent By taking the quotients of $\RR_G^x$ we get the
following 3-Sasakian, twistor and QK reductions:
$$\Z_3\setminus G_2(V)/{\rm Sp}(1),$$
$$\Z_3\setminus G_2(V)/{\rm Sp}(1){\rm U}(1),$$
$$\Z_3\setminus G_2(V)/{\rm SO}(4).$$
In particular, the QK reduction in the compact case as been
already found in \cite{KoSw93}. So, the reduced diamond diagram is
covered (with branches) by the associated diamond diagram to the
Wolf space $G_2(V)/{\rm SO}(4)$. Note that the action of $T_x$ at
the 3-Sasakian level is not locally free, but only
\emph{quasi-free}. In fact, this is a typical case of a
not-locally free action which produces anyway an orbifold
structure on the reduction, thanks to existence of a
\emph{section}:
$$
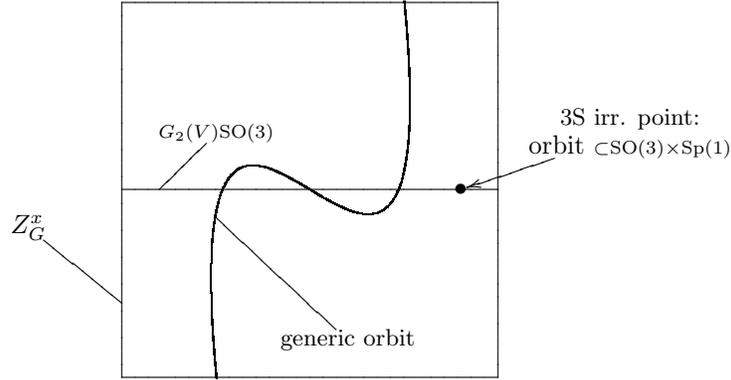
\xy
<12cm,-2.5cm>;<12.25cm,-2.5cm>:(-10,-10)="A";(10,-10)="B";(10,10)="C";(-10,10)="D";%
"A";"B"**{\dir{-}};"C"**{\dir{-}};"D"**{\dir{-}};"A"**{\dir{-}};%
"A";"B"**{}?<(.25)="ab"; "C";"D"**{}?<(.25)="cd";%
"B";"C"**{}?="m1";"A";"D"**{}?="m2"; "m1";"m2"**{\dir{-}};%
"ab";"cd"**\crv{(-7.5,17.5)&(7.5,-17.5)}?(.15)="z";%
"m1";"m2"**{}?<(.1)*{\bullet}="s";%
(17,3)*{\scriptstyle \txt{\small 3S irr. point:\\ orbit $\scriptstyle\subset {\rm SO}(3)\times {\rm Sp}(1)$}}="p";%
(-15,-2)*{Z_G^x}="q";%
(0,5)="t";%
"A";"D"**{}?>(.2)="u";%
"m1";"m2"**{}?<(.9)="v";%
(-5,3)*{\scriptstyle G_2(V){\rm SO}(3)}="w";%
(2,-8)*{ \txt{\small generic orbit}}="x";
\ar @{->}"p";"s"%
\ar @{-}"q";"u"%
\ar @{-}"w";"v"%
\ar @{-}"x";"z"%
\endxy$$

\noindent The branched locus of the reduction consists, at the
level of the group, of the set of all the (3-Sasakian) irregular
points in ${\rm SO}(V)$.\\
In order to describe it, we can consider its quotient by ${\rm
SO}(3)$, which is a subset of the Stiefel manifold
${\VV}_{0,4}(V)$.
\begin{teo}Let us define $z_r:=e_{2r}+ie_{2r+1}\in V\otimes\C$.
Then, the quotient of the fixed points set by  ${\rm SO}(3)$ is
$$({\rm U}(z_1,z_2,z_3)\cdot{\rm Sp}(1)_-)/{\rm U}(z_1)\subset {\VV}_{0,4}(V)$$
Hence, the 3-Sasakian branch locus is
$${\rm U}(z_1,z_2,z_3)/({\rm U}(z_1)\times \Delta({\rm U}(z_2,z_3))),$$
where $\Delta({\rm U}(z_2,z_3)))$ denotes the subgroup of ${\rm
U}(z_2,z_3)$ which consists of the complex multiples of the
identity. Finally, the quaternionic branch locus is
$${\rm U}(z_1,z_2,z_3)/({\rm U}(z_1)\times {\rm U}(z_2,z_3)).$$
\end{teo}
\noindent Now, let us consider the noncompact case with $x\in
S^+\cup K-\{0\}$. In this case, the picture is completely
different! In fact, the set $G_{2(2)}{\rm SO}(3)$ isn't a section
anymore, but the action is proper and free, so we obtain
(nonempty) \emph{manifolds}. As we'll see, we get whole families
of \emph{manifolds} carrying out the reduction.

\subsection{The adjoint orbits of ${\rm SO}(3,4)$ and the associated reductions}

\noindent In the following we specialize a general construction in
\cite{BurCush} to the case of orthogonal groups. An analogous
application of the construction (in the case of symplectic groups)
can be found in \cite{BCGP05}.\\

\noindent For any complex vector space $V$ and a symmetric,
complex-valued, non-degenerate bilinear form $\tau$ let ${\rm
O}(V,\tau)$ be the group of all the (linear, complex)
automorphisms of $V$ which preserve the bilinear form $\tau$ and
let $\mathfrak{o}(V,\tau)$ be its Lie algebra, which is the
algebra of all skew-selfadjoint complex endomorphisms of $V$.\\
Let $\sigma:V\rightarrow V$ be an anti-endomorphism (i.e.,
$\sigma(\alpha v)=\overline{\alpha}\sigma(v)$ for any $\alpha\in
\C$, $v\in V$) such that $\sigma^2=1$,
$\tau^{\sigma}=\overline{\tau}$.\\
If $V_{\sigma}:=\{v\in V:\sigma(v)=v\}$ and $\tau_{\sigma}$ is the
restriction of $\tau$ on this real subspace, then $\tau_{\sigma}$
is real and non-degenerate. The group
$${\rm O}(V,\tau,\sigma):=\{g\in {\rm O}(V,\tau):g^{\sigma}=g\}$$
can be identified (trough restriction to $V_{\sigma}$) with the
group of real automorphisms of $V_{\sigma}$ which preserve
$\tau_{\sigma}$, and so it is isomorphic to ${\rm O}(k,l)$, where
$(k,l)$ is the signature of $\tau_{\sigma}$. Let
${\mathfrak{o}}(V,\tau,\sigma)$ be its Lie algebra.\\
Let $A\in {\mathfrak{o}}(V,\tau,[\sigma])$, $A'\in
{\mathfrak{o}}(V',\tau',[\sigma'])$. We say $A,A'$ to be
\emph{equivalent} if there exists an isomorphism
$\phi:V\rightarrow V'$ such that $\tau'=\tau^{\phi}$,
[$\sigma'=\phi\sigma\phi^{-1}$] and
$A'=\phi A\phi^{-1}$.\\
This definition is actually an equivalence relation and its
classes are called \emph{types}:$$\Delta:=[A].$$ Furthermore,
$\phi$ as above defines an isomorphism between ${\rm
O}(V,\tau,[\sigma])$ and ${\rm O}(V',\tau',[\sigma'])$, and in the
case $(V,\tau,[\sigma])=(V',\tau',[\sigma'])$, $A,A'$ are
equivalent if and only if they belong to the same adjoint orbit.
It's straightforward to define the \emph{sum} of two types:
$$\Delta,\Delta'\mapsto \Delta\oplus \Delta'.$$
\begin{teo}Every type can be written (up to the order) in an unique way as a sum of \emph{indecomposable} types.\end{teo}
Every $A\in{\mathfrak{o}}(V,\tau,[\sigma])$ can be decomposed in
an unique way as a sum
$$A=S+N,$$
where $S,N\in{\mathfrak{o}}(V,\tau,[\sigma])$, $[S,N]=0$ and are
respectively semisimple and nilpotent. The nonnegative integer $k$
such that $N^k\neq 0$, $N^{k+1}=0$ is called \emph{the height} of
$A$. It is an invariant of the type $\Delta=[A]$ and we denote it
with ${\rm ht}(\Delta)$. (Another invariant is the
\emph{dimension} of the type.)
A type with null height is called \emph{semisimple}.\\
There are two kinds of semisimple indecomposable types for ${\rm
O}(V,\tau)$:\begin{itemize} \item ${\Delta}(\zeta,-\zeta)$,
$\zeta\neq 0$: it is  $\C$-generated by two nonorthogonal vectors
$v,w$ of norm $0$, eigenvectors of $S$ with respect to the
eigenvectors $\zeta,-\zeta$; \item ${\Delta}(0)$, if $S=0$: it is
$\C$-generated by a vector of norm $1$.
\end{itemize}
In the case of ${\rm O}(V,\tau,\sigma)$ we get more semisimple
indecomposable types, namely
\begin{itemize}
\item
$\Delta(\zeta,-\zeta,\overline{\zeta},-\overline{\zeta})$\quad
with \quad $\overline{\zeta}\neq\pm\zeta$. This \emph{real} type
can be construct as follows: take the complex type
$$\Delta(\zeta,-\zeta)\oplus\Delta(\overline{\zeta},-\overline{\zeta})$$
which is $\C$-generated by $v_1,v_2,v_3,v_4$ such that
$$\tau(v_h,v_l)=\left\{\begin{array}{lc}1&\mbox{ if }\ \ \{h,l\}=\{1,2\}\mbox{ or }\{3,4\};\\
                                     &                                             \\
                                    0&\mbox{ otherwise }\end{array}\right.$$
and which are eigenvectors of
$\zeta,-\zeta,\overline{\zeta},-\overline{\zeta}$, respectively.
Furthermore, after rescaling we can suppose $V_{\sigma}$ is
$\R$-generated by the orthogonal vectors
$$w_1=\frac{1}{2}(v_1+v_2+v_3+v_4),\quad\|w_1\|=1,\ $$
$$w_2=\frac{i}{2}(v_1-v_2-v_3+v_4),\quad\|w_2\|=1,\ $$
$$w_3=\frac{1}{2}(v_1-v_2+v_3-v_4),\quad\|w_3\|=-1,$$
$$w_4=\frac{i}{2}(v_1+v_2-v_3-v_4),\quad\|w_4\|=-1.$$
Hence, we get
$$\left\{\begin{array}{lcl}Sw_1&=&aw_3+bw_2\\
                           Sw_2&=&aw_4-bw_1\\
                           Sw_3&=&aw_1+bw_4\\
                           Sw_4&=&aw_2-bw_3,\end{array}\right.$$
where $\zeta=a+ib$, $a,b\in \R$. \item $\Delta(\zeta,-\zeta)$\quad
with\quad $\overline{\zeta}=\zeta\neq 0$: Let $v_1,v_2$ be the
generators of the respective complex type. In this case
$V_{\sigma}$ is $\R$-generated by $v_1,v_2$ and also by the
orthogonal vectors
$$w_1:=\frac{1}{\sqrt{2}}(v_1+v_2),\quad\|w_1\|=1,\ $$
$$w_2:=\frac{1}{\sqrt{2}}(v_1-v_2),\quad\|w_2\|=-1.$$
Hence, we get
$$\left\{\begin{array}{lcl}Sw_1&=&aw_2\\
                           Sw_2&=&aw_1,\end{array}\right.$$
where $\zeta=a+i0$, $a\in \R$.

\item  $\Delta^{\pm}(\zeta,-\zeta)$\quad with
\quad$\overline{\zeta}=-\zeta\neq 0$: As in the previous case,
$V_{\sigma}$ is generated by two orthogonal vectors $w_1,w_2$ such
that
$$\|w_1\|=\|w_2\|=\pm 1$$
and
$$\left\{\begin{array}{lcl}Sw_1&=&bw_2\\
                           Sw_2&=&-bw_1,\end{array}\right.$$
where $\zeta=ib$, $b\in \R$. \item ${\Delta}^{\pm}(0)$ with $S=0$:
In this case, $V_{\sigma}$ is generated by a vector of norm $\pm
1$.
\end{itemize}
Let $\Delta$ be a generic indecomposable type, of height $k$. Then
${\rm Im}\ N\leq{\rm Ker} \ N^k$, while we say $\Delta$ to be
\emph{uniform} if the equality holds. If $\Delta$ is uniform, we
define
$$\widetilde{V}:=V/{\rm Im}\ N,$$
$$\widetilde{A}(v+{\rm Im}\ N):=Av+{\rm Im}\ N=Sv+{\rm Im}\ N,$$
$$\widetilde{\sigma}(v+{\rm Im}\ N)=\sigma(v)+{\rm Im}\ N,$$
$$\widetilde{\tau}(v+{\rm Im}\ N,w+{\rm Im}\ N):=\tau(v,N^k(w)).$$
In particular, the bilinear form $\widetilde{\tau}$ is
non-degenerate and it is symmetric if the height is even: in this
case, we get so a semisimple type $\widetilde{\Delta}$. In the odd
case we no longer get an orthogonal type but a \emph{symplectic}
one (i.e., substitute in the previous definitions ${\rm
O},\mathfrak{o}$,\emph{symmetric} with ${\rm Sp},\mathfrak{sp}$,
\emph{skewsymmetric}). However, \begin{enumerate} \item if
$\Delta$ is indecomposable, then it is uniform and
$\widetilde{\Delta}$ is indecomposable; \item if $\Delta$ is
uniform, then it is uniquely determined by its height and
$\widetilde{\Delta}$.
\end{enumerate}
In addition we need to describe a classification of the
indecomposable semisimple types for ${\rm Sp}(V,\tau,[\sigma])$.
\\All the indecomposable semisimple types of ${\rm Sp}(V,\tau)$
can be denoted again with ${\Delta}(\zeta,-\zeta)$, where
$\zeta\in \C$: in this type, $V$ is $\C$-generated by two vectors
$v,w$ which are eigenvectors of $\zeta,-\zeta$ respectively.
Moreover, the  indecomposable semisimple types of ${\rm
Sp}(V,\tau,\sigma)$ are\begin{itemize} \item
$\Delta(\zeta,-\zeta,\overline{\zeta},-\overline{\zeta})$\quad
with \quad $\overline{\zeta}\neq\pm\zeta$: As before, $V_{\sigma}$
is $\R$-generated by $w_1,w_2,w_3,w_4$ such that
$$\left\{\begin{array}{lcl}Sw_1&=&aw_3+bw_2\\
                           Sw_2&=&aw_4-bw_1\\
                           Sw_3&=&aw_1+bw_4\\
                           Sw_4&=&aw_2-bw_3,\end{array}\right.$$
where $\zeta=a+ib$, $a,b\in \R$. Furthermore, $\tau(w_h,w_l)=-1$
if $(h,l)=(1,2)$ or $(3,4)$, and it is equal to $0$ otherwise.
\item $\Delta(\zeta,-\zeta)$\quad with\quad
$\overline{\zeta}=\zeta\neq 0$: Here $V_{\sigma}$ is
$\R$-generated by $w_1,w_2$ such that
$$\left\{\begin{array}{lcl}Sw_1&=&aw_2\\
                           Sw_2&=&aw_1,\end{array}\right.$$
where $\zeta=a+i0$, $a\in \R$, and $\tau(w_1,w_2)=-1$. \item
$\Delta^{\pm}(\zeta,-\zeta)$\quad
with\quad$\overline{\zeta}=-\zeta\neq 0$: In this case,
$V_{\sigma}$ is $\R$-generated by $w_1,w_2$ such that
$$\left\{\begin{array}{lcl}Sw_1&=&bw_2\\
                           Sw_2&=&-bw_1,\end{array}\right.$$
where $\zeta=ib$, and $\tau(w_1,w_2)=\pm 1$. \item $\Delta(0,0)$
with $S=0$: In this case, $V_{\sigma}$ has real dimension 2.
\end{itemize}
It remains to describe how $\Delta$ can be recovered from
$\widetilde{\Delta}$ and its height. Let $\Delta$ be an
indecomposable ($\Rightarrow$ uniform) type of height $k$. Then
$$V=W\oplus NW\oplus N^2W\oplus\cdots\oplus N^kW$$
such that the subspace $W$ is $S$-invariant, $\sigma$-invariant
and $W\bot N^hW$ for $h=0,\cdots, k-1$. Furthermore, the subspaces
$N^hW$ are isomorphic and the restriction of $S$ to $W$ (endowed
with $\widetilde{\tau}, \widetilde{\sigma}$) represents the type
$\widetilde{\Delta}$.\\
Conversely, let $\Delta=[S]$ a semisimple indecomposable type and
$k$ a positive integer. If $S\in \gg(V,\tau,\sigma)$, let us
define
$$V_k:=\bigoplus^kV,\qquad
S_k:=\bigoplus^kS,\qquad \sigma_k:=\bigoplus^k\sigma.$$
Furthermore, let $N_k$ be the endomorphism of $V_k$ which does
shift the components, namely
$$N_k(v_0,v_1,\cdots,v_k):=(0,v_0,\cdots,v_{k-1})$$
and let be $\tau_k$ the unique bilinear form on $V_k$ with respect
to which $N$ is skewsymmetric and such that, for every $u,v\in W$,
$$\tau_k(u,N^iv)=\left\{\begin{array}{lc}\tau(u,v)&\mbox{ if } i=k\\&\\
                                           0&\mbox{ otherwise. }\end{array}\right.$$
Then, $\Delta_k:=[S_k+N_k]$ is an indecomposable type of height
$k$ such that $\widetilde{\Delta_k}=\Delta$. Below is the list of
all the indecomposable, orthogonal types of dimension $\leq 7$:
\begin{table}[h] \centering
\begin{tabular}{|c|c|c|c|}\hline $\Delta$ & & ${\rm dim}({\Delta})$& signature\\
                          \hline $\Delta_6^-(0)$ & & $7$& $(4,3)$\\
                          \hline $\Delta_6^+(0)$ & & $7$& $(3,4)$\\
                          \hline $\Delta_2^-(\zeta,-\zeta)$ &$0\neq\zeta\in i\R$ & $6$& $(4,2)$\\
                          \hline $\Delta_2(\zeta,-\zeta)$ & $0\neq\zeta\in\R$ & $6$& $(3,3)$\\
                          \hline $\Delta_2^+(\zeta,-\zeta)$ &$0\neq\zeta\in i\R$ & $6$& $(2,4)$\\
                          \hline $\Delta_4^+(0)$ & & $5$& $(3,2)$\\
                          \hline $\Delta_4^-(0)$ & & $5$& $(2,3)$\\
                          \hline $\Delta_0(\zeta,-\zeta,\overline{\zeta},-\overline{\zeta})$ &$\zeta\in \C-(\R\cup i\R)$ & $4$& $(2,2)$\\
                          \hline $\Delta_1^-(\zeta,-\zeta)$ &$0\neq\zeta\in i\R$ & $4$& $(2,2)$\\
                          \hline $\Delta_1(\zeta,-\zeta)$ & $\zeta\in\R$ & $4$& $(2,2)$\\
                          \hline $\Delta_1^+(\zeta,-\zeta)$ &$0\neq\zeta\in i\R$ & $4$& $(2,2)$\\
                          \hline $\Delta_2^-(0)$ & & $3$& $(2,1)$\\
                          \hline $\Delta_2^+(0)$ & & $3$& $(1,2)$\\
                          \hline $\Delta_0^+(\zeta,-\zeta)$ &$0\neq\zeta\in i\R$ & $2$& $(2,0)$\\
                          \hline $\Delta_0(\zeta,-\zeta)$ & $0\neq\zeta\in\R$ & $2$& $(1,1)$\\
                          \hline $\Delta_0^-(\zeta,-\zeta)$ &$0\neq\zeta\in i\R$ & $2$& $(0,2)$\\
                          \hline $\Delta_0^+(0)$ & & $1$& $(1,0)$\\
                          \hline $\Delta_0^-(0)$ & & $1$& $(0,1)$\\
\hline
                \end{tabular} \end{table}
\ \\By combining them together, we get a list of 48 orthogonal
types with signature $(3,4)$. \\Furthermore, by setting
\begin{eqnarray}\Delta_2^+(0,0)&:=&2\Delta_2^+(0)\nonumber\\
                \Delta_2(0,0)&:=&\Delta_2^+(0)\oplus\Delta_2^-(0)\nonumber\\
                \Delta_0(a,-a,a,-a)&:=&2\Delta_0(a,-a)\nonumber\\
                \Delta_0(ib,-ib,-ib,ib)&:=&\Delta_0^+(ib,-ib)\oplus\Delta_0^-(ib,-ib)\nonumber\\
                \Delta_0(0,0)&:=&\Delta_0^+(0)\oplus\Delta_0^-(0)\nonumber\\
                \Delta_0^{\pm}(0,0)&:=&2\Delta_0^{\pm}(0)\nonumber\end{eqnarray}
for $a,b\in\R$, we obtain only 24 families:
\begin{table}[h] \centering
\begin{tabular}{|c|c|c|c|}\hline $\Delta$ & ${\rm ht}(\Delta)$& $\#$ paramts. & Name\\
                          \hline $\Delta_6^+(0)$& 6 & 0 & ${\rm I}_1$\\
                          \hline $\Delta_4^+(0)\oplus\Delta_0^-(\zeta,-\zeta)$& 4 & 1 & ${\rm II}_1$\\
                          \hline $\Delta_4^-(0)\oplus\Delta_0(\zeta,-\zeta)$& 4 & 1 & ${\rm II}_2$\\
                          \hline $\Delta_2(\zeta,-\zeta)\oplus\Delta_0^-(0)$& 2 & 1 & ${\rm II}_3$\\
                          \hline $\Delta_2^+(\zeta,-\zeta)\oplus\Delta_0^+(0)$& 2 & 1 & ${\rm II}_4$\\
                          \hline $\Delta_1^-(\zeta,-\zeta)\oplus\Delta_2^+(0)$& 2 & 1 & ${\rm II}_5$\\
                          \hline $\Delta_1(\zeta,-\zeta)\oplus\Delta_2^+(0)$& 2 & 1 & ${\rm II}_6$\\
                          \hline $\Delta_1^+(\zeta,-\zeta)\oplus\Delta_2^+(0)$& 2 & 1 & ${\rm II}_7$\\
                          \hline $\Delta_0(\zeta,-\zeta,\overline{\zeta},-\overline{\zeta})\oplus\Delta_2^+(0)$& 2 & 2 & ${\rm III}_1$\\
                          \hline $\Delta_2^-(0)\oplus\Delta_0(\zeta_1,-\zeta_1)\oplus\Delta_0^-(\zeta_2,-\zeta_2)$& 2 & 2 & ${\rm III}_2$\\
                          \hline $\Delta_0^+(\zeta_1,-\zeta_1)\oplus\Delta_2^+(0)\oplus\Delta_0^-(\zeta_2,-\zeta_2)$& 2 & 2 & ${\rm III}_3$\\
                          \hline $\Delta_2^+(0)\oplus \Delta_0(\zeta_1,-\zeta_1)\oplus\Delta_0(\zeta_2,-\zeta_2)$& 2 & 2 & ${\rm III}_4$\\
                          \hline $\Delta_1^-(\zeta_1,-\zeta_1)\oplus\Delta_0(\zeta_2,-\zeta_2)\oplus\Delta_0^-(0)$& 1 & 2 & ${\rm III}_5$\\
                          \hline $\Delta_1(\zeta_1,-\zeta_1)\oplus\Delta_0(\zeta_2,-\zeta_2)\oplus\Delta_0^-(0)$& 1 & 2 & ${\rm III}_6$\\
                          \hline $\Delta_1^+(\zeta_1,-\zeta_1)\oplus\Delta_0(\zeta_2,-\zeta_2)\oplus\Delta_0^-(0)$& 1 & 2 & ${\rm III}_7$\\
                          \hline $\Delta_1^-(\zeta_1,-\zeta_1)\oplus\Delta_0^+(0)\oplus\Delta_0^-(\zeta_2,-\zeta_2)$& 1 & 2 & ${\rm III}_8$\\
                          \hline $\Delta_1(\zeta_1,-\zeta_1)\oplus\Delta_0^+(0)\oplus\Delta_0^-(\zeta_2,-\zeta_2)$& 1 & 2 & ${\rm III}_9$\\
                          \hline $\Delta_1^+(\zeta_1,-\zeta_1)\oplus\Delta_0^+(0)\oplus\Delta_0^-(\zeta_2,-\zeta_2)$& 1 & 2 & ${\rm III}_{10}$\\
                          \hline $\Delta_0(\zeta_1,-\zeta_1,\overline{\zeta_1},-\overline{\zeta_1})\oplus\Delta_0(\zeta_2,-\zeta_2)\oplus\Delta_0^-(0)$& 0 & 3 & ${\rm IV}_1$\\
                          \hline $\Delta_0(\zeta_1,-\zeta_1,\overline{\zeta_1},-\overline{\zeta_1})\oplus\Delta_0^+(0)\oplus\Delta_0^-(\zeta_2,-\zeta_2)$& 0 & 3 & ${\rm IV}_2$\\
                          \hline $\Delta_0^+(\zeta_1,-\zeta_1)\oplus\Delta_0(\zeta_2,-\zeta_2)\oplus\Delta^-_0(\zeta_3,-\zeta_3)\oplus \Delta_0^-(0)$& 0 & 3 & ${\rm IV}_3$\\
                          \hline $\Delta_0^+(\zeta_1,-\zeta_1)\oplus\Delta_0^+(0)\oplus \Delta^-_0(\zeta_2,-\zeta_2)\oplus\Delta^-_0(\zeta_3,-\zeta_3)$& 0 & 3 & ${\rm IV}_4$\\
                          \hline $\Delta_0(\zeta_1,-\zeta_1)\oplus\Delta_0(\zeta_2,-\zeta_2)\oplus\Delta_0(\zeta_3,-\zeta_3)\oplus\Delta_0^-(0)$& 0 & 3 & ${\rm IV}_5$\\
                          \hline $\Delta_0(\zeta_1,-\zeta_1)\oplus\Delta_0(\zeta_2,-\zeta_2)\oplus\Delta_0^+(0)\oplus\Delta_0^-(\zeta_3,-\zeta_3)$& 0 & 3 & ${\rm IV}_6$\\
\hline
                \end{tabular} \end{table}

In particular,
$$A_x\in\left\{\begin{array}{lc}\Delta_0^{+}(i,-i)\oplus 2\Delta_0^{-}(-i,i)\oplus\Delta_0^{+}(0)&\mbox{ if
}\|x\|=1,\\&\\
                                3\Delta_0(1,-1)\oplus\Delta_0^{-}(0)&\mbox{ if
}\|x\|=-1,\\&\\
                                \Delta_1(0,0)\oplus\Delta_2^+(0)&\mbox{ if }\|x\|=0.\end{array}\right.$$

One can easily show
\begin{teo}Let $A\in \mathfrak{so}(3,4)$, belonging to the adjoint orbit $[A]$. If $[A]\notin {\rm IV}_4$, then the quaternionic
action generated by $A$ is proper and free at all levels.\end{teo}
\emph{Proof:} If $[A]\notin {\rm IV}_4$,  the matrix ${\rm
exp}(tA)\in {\rm SO}(3,4)$ contains polynomial expressions and/or
hyperbolic functions, so
$$\lim_{t\rightarrow\pm\infty}\|{\rm exp}(tA)\|=+\infty,$$
where $\|\cdot\|$ denotes the euclidean norm of $\R^{49}$. So, the
group generated by $A$ is the Lie subgroup of ${\rm SO}(3,4)$ and
the action is proper on the group level. Furthermore, the
restriction of the action on $Z_G$ is proper, as well as the
actions induced on the 3-Sasakian, the twistor, and the QK level.\\
Since the group generated by $A$ is isomorphic to $\R$, the
freeness follows from the properness.\qed
\begin{cor}Let $A\in \mathfrak{so}(3,4)$, belonging to the adjoint orbit $[A]$. If $[A]\notin {\rm IV}_4$, then the reduction associated to
$A$ produces \emph{manifolds}. \end{cor}\noindent Furthermore,
from a straightforward computation of the centralizers, it follows
that the reduced diamond diagrams admit (at least) two commuting
Killing vector fields.

\begin{teo}Let $A\in \mathfrak{so}(3,4)$, belonging to the adjoint orbit $[A]={\rm IV}_4(a,b,c)$. The quaternionic
action generated by $A$ is proper if and only if the triple
$(a,b,c)$ is commensurable. In this case, if $0\notin\{a,b,c\}$ we
can suppose $a,b,c$ to be positive integers with ${\rm
gcd}(a,b,c)=1$. Furthermore, there are no 3-Sasakian irregular
points if and only if $a,b,c$ are distinct.\end{teo}

\noindent Analogously to the compact case, it is possible to give
conditions to the coefficients $(a,b,c)$ in order to obtain a free
3-sasakian action (see \cite{BGP02}). It seems hard to can find
stronger conditions for the QK action.

\begin{rem}
The previous results can be used to study the 1-dimensional QK
reductions in the case of the Wolf space $G_{2(2)}/{\rm SO}(4)$.
In fact, any element $v\in \gg_{2(2)}$ determines an adjoint orbit
in $\mathfrak{so}(3,4)$ $\Delta$. So, $\Delta\notin {\rm IV}_4$ if
and only if the adjoint orbit of $v$ in $\gg_{2(2)}$ doesn't
intersect $\mathfrak{so}(4)\subset \gg_{2(2)}$ and in this case
the reduction produces self-dual, Einstein \emph{4}-manifolds with
negative sectional curvature and \emph{few} symmetries.

\end{rem}
\noindent\textbf{Open problems:}\begin{itemize} \item In the
non-${\rm IV}_4$ cases, classify the obtained manifolds. In
particular, is it true that the quaternionic Killing algebra is
the quotient $\mathfrak{cent}(\gt)/\gt$? Does it coincide with the
Killing algebra of the whole reduced diamond diagram? Are there
algebraic invariants which classify this manifolds?

\item In the ${\rm IV}_4$ case, classify the obtained orbifolds.
In the case of two equal parameters, does the action admit
sections (like in the case $(1,1,1)$)?

\item Do the same for 2-dimensional toric actions. How to classify
the conjugacy classes of abelian subalgebras (of dimension $>1$)
of classic Lie algebras?

\item Repeat the whole construction for the other Wolf spaces. In
particular, we want to find $4$-dimensional manifolds with few
symmetries, so would be better to start from low-dimensional,
noncompact Wolf spaces, for example $G_{2(2)}/{\rm SO}(4)$. How to
classify the adjoint orbits of $G_{2(2)}$ and in general, of an
exceptional Lie group?

\item In the case of non-Riemannian Wolf spaces, reduction process
can encounter a further problem: in fact, the \emph{reduced metric}
may have some kernel.

\end{itemize}
\date{}

\bibliographystyle{plain}
\bibliography{wolf}
\end{document}